\documentclass[11pt]{article}
\usepackage{epsfig}
\usepackage{amssymb}
\usepackage{amsmath}
\newtheorem{theorem}{Theorem}[section]

\newtheorem{proposition}[theorem]{Proposition}
\newtheorem{corollary}[theorem]{Corollary}

\newtheorem{example}[theorem]{Example}
\newtheorem{remark}[theorem]{Remark}
\numberwithin{equation}{section}
\def\qed{\hfill $\Box$\medskip}
\def\IR{{\bf R}}
\def\IC{{\bf C}}

\def\bA{{\mathbf {A}}}
\def\bB{{\bf B}}

\def\bF{{\bf F}}

\def\bV{{\bf V}}
\def\ba{{\bf a}}
\def\bb{{\bf b}}
\def\bc{{\bf c}}

\def\b0{{\bf 0}}
\def\bmu{{\pmb \mu}}

\def\v{{\bf v}}
\def\x{{\bf x}}
\def\y{{\bf y}}
\def\z{{\bf z}}
\def\cX{{\cal X}}
\def\cY{{\cal Y}}

\def\cS{{\cal S}}

\def\cB{{\cal B}}

\def\cU{{\cal U}}
\def\cF{{\cal F}}

\def\cP{{\cal P}}
\def\cV{{\cal V}}
\def\cH{{\cal H}}
\def\BH{{\cB(\cH)}}
\def\FH{{\cF(\cH)}}
\def\SH{{\cS(\cH)}}
\def\cl{{\bf cl}\,}

\def\conv{{\bf conv}\,}

\def\diag{{\rm diag}\,}

\def\rank{{\rm rank}\,}
\def\span{{\rm span}\,}
\def\range{{\rm range}\,}
\def\dim{{\rm dim}\,}
\def\ker{{\rm ker}\,}
\def\la{{\langle}}
\def\ra{{\rangle}}
\def\Ra{\Rightarrow }

\def\({\left(}
\def\){\right)}
\def\[{\left[}
\def\]{\right]}
\def\dfrac{\displaystyle\frac}

\begin{document}

\title{Quantum error correction and generalized numerical ranges}
\author{{\sc Chi-Kwong Li}
and
{\sc Yiu-Tung Poon}}
\date{}

\maketitle

\medskip\noindent{\large{\bf Abstract}}

For a noisy quantum channel,
a quantum error correcting code exists
if and only if the joint higher rank
numerical ranges associated with the error operators of the
channel is non-empty. In this paper, geometric properties of
the joint higher rank numerical ranges are obtained and their
implications to quantum computing are discussed.
It is shown that if the dimension of the underlying
Hilbert space of the quantum states is sufficiently large,
the joint higher rank numerical range
of operators  is always
star-shaped and contains a non-empty convex subset.
In case the operators are infinite
dimensional, the joint  infinite rank numerical range
of the operators is a convex set 
closely related to the joint essential numerical ranges of
the operators. 

\medskip
\bf AMS Subject Classification \rm 47A12, 15A60, 15A90, 81P68.

\bf Keywords \rm Quantum error correction,
joint higher rank numerical range, joint essential numerical range,
self-adjoint operator, Hilbert space.

\section{Introduction}

In quantum computing, information is stored in
{\it quantum bits}, abbreviated as {\it qubits}.
Mathematically, a qubit is represented by a
$2\times 2$ rank one Hermitian matrix $Q = vv^*$, where
$v \in \IC^2$ is a unit vector.
A state of $N$-qubits $Q_1, \dots, Q_N$ is represented by their
tensor products in $M_n$ with $n = 2^N$.  A {\it quantum channel}
for states of
$N$-qubits corresponds to a {\it trace preserving completely positive linear
map} $\Phi: M_n \rightarrow M_n$. By the structure theory of
completely positive linear map \cite{C},
there are $T_1, \dots, T_r \in M_n$
with $\sum_{j=1}^r   T_j^*T_j = I_n$ such that
\begin{equation}\label{channel}
\Phi(X) = \sum_{j=1}^r T_jXT_j^*.
\end{equation}
In the context of quantum error correction, $T_1, \dots, T_r$ are known
as the {\it error operators}.

Let $\bV$ be a $k$-dimensional subspace of $\IC^n$ and $P$
the orthogonal projection of
$\IC^n$ onto $\bV$. Then $\bV$
is a {\it quantum error correcting code} for the quantum channel
$\Phi$ if there exists another trace
preserving completely positive linear map $ \Psi: M_n
\rightarrow M_n$ such that  $\Psi\circ\Phi(A)=A$ for all $A\in PM_nP$.
By the results in \cite{KL},
this happens if and only if there are scalars $\gamma_{ij}$
with $1 \le i, j \le r$  such that
$$PT_i^*T_jP = \gamma_{ij}P.$$
Let $\cP_k$ be the set of rank $k$ orthogonal projections
in $M_n$. Define the {\it joint rank $k$-numerical range}
of an $m$-tuple of matrices $\bA = (A_1, \dots, A_m) \in M_n^m$  by
$$\Lambda_k(\bA) = \{(a_1, \dots, a_m) \in \IC^m:
\hbox{ there is } P \in \cP_k \hbox{ such that }
PA_jP = a_j P \hbox{ for } j = 1, \dots, m \}.$$
Then the quantum channel $\Phi$ defined in (\ref{channel})
has an error correcting code of $k$-dimension
if and only if
$$\Lambda_k(T_1^*T_1, T_1^*T_2, \dots, T_r^*T_r) \ne \emptyset.$$
Evidently,
$(a_1, \dots, a_m)\in \Lambda_k(\bA)$
if and only if there exists an $n\times k$
matrix $U$ such that
$$U^*U=I_k, \quad \hbox{  and  } \quad
U^*A_jU=a_jI_k  \qquad \hbox{ for } j = 1, \dots, m.$$
Let $\x,\y\in\IC^n$. Denote by $\la \bA\x,\y\ra $
the vector $\(\la A_1\x,\y\ra ,\dots,\la A_m\x,\y\ra \)\in\IR^m$.
Then $
\ba\in \Lambda_k(\bA)$ if and only if there exists
an orthonormal set $\{\x_1,\dots,\x_k\}$ in $\IC^n$
such that
$$\la \bA \x_i,\x_j\ra=\delta_{ij}\ba,$$
where $\delta_{ij}$ is the Kronecker delta.
When $k = 1$, $\Lambda_1(\bA)$ reduces to the
(classical)
{\it joint numerical range}
$$W(\bA) = \{(\x^*A_1\x, \dots, \x^*A_m\x):\x \in \IC^n, \ \x^*\x = 1\}$$
of $\bA$,  which is quite
well studied; see \cite{LP} and the references therein.
It turns out that even for a single matrix
$A \in M_n$, the study of $\Lambda_k(A)$ is highly non-trivial,
and the results are useful in quantum computing, say,
in constructing binary unitary channels;
see \cite{Cet,Cet0,Cet1,Cet2,GLW,LPS,LS,W1}.

More generally, let $\BH$ be the algebra of bounded linear operators
    acting on a Hilbert space $\cH$, which may be infinite dimensional.
    One can extend the definition of $\Lambda_k(A)$ to $A \in \BH$.
    If $\cH$ is infinite dimensional, one may allow $k = \infty$ by
    letting $\cP_k$ be the set of infinite rank orthogonal projections
    in $\BH$ in the definition; see \cite{LPS2,R}. There are a number of reasons to consider rank $k$-numerical range
of infinite dimensional operators. First, many quantum
mechanical phenomena are better described using infinite
dimensional Hilbert spaces. Also, a practical quantum computer
must be able to handle a large number of qubits so that the
underlying Hilbert space must have a very large dimension.
We will also consider the joint rank
$k$-numerical range of an $m$-tuple $\bA = (A_1, \dots, A_m)$
of infinite dimensional
operators $A_1, \dots, A_m$ for positive integers $k$ and $k = \infty$.
It is interesting to note that
$\Lambda_\infty(\bA)$ has intimate connection with the
{\it joint essential numerical range} of $\bA$ defined as
$$W_e(\bA) = \cap \{ \cl(W(\bA + \bF)): \bF \in \FH^m\},$$
where $\FH$ denotes the set of finite rank operators in $\BH$
and $\cl(\cS)$ denotes the closure of the set $\cS$.
Clearly, the joint essential numerical range
is useful for the study of the joint behaviors of operators under
perturbations of finite rank (or compact) operators.

The purpose of this paper is to study the joint rank $k$-numerical range
of $\bA = (A_1, \dots, A_m) \in \BH^m$.
Understanding the properties of $\Lambda_k(\bA)$ is useful for constructing
quantum error correcting codes and studying their properties
such as their stability under perturbation.

Our paper is organized as follows.
In Section 2, we present some basic properties of
$\Lambda_k(\bA)$.
Section 3 concerns the geometric properties
of $\Lambda_k(\bA)$.
We show that if $\dim \cH$ is sufficiently large,
then $\Lambda_k(\bA)$ is always star-shaped and contains a
convex subset.
In Section 4, we study the connection between  $W_e(\bA)$,
$\Lambda_k(\bA)$ and its closure $\cl(\Lambda_k(\bA))$.
We show that $\Lambda_\infty(\bA)$ is always convex, and is a subset
of the set of star centers of $\Lambda_k(\bA)$ for each positive
integer $k$. We also show that
$$
W_e(\bA) = \cap_{k\ge 1} \cl(\Lambda_k(\bA))\,.$$
Moreover, we obtain several equivalent formulations of
$\Lambda_\infty(\bA)$ including
$$\Lambda_\infty(\bA)
= \cap \{\Lambda_k(\bA+\bF): \bF \in \FH^m\}.$$
The results extend those in \cite{AS,LP2}.

Let $\SH$ be the real linear space of self-adjoint operators
in $\BH$.
Suppose
$$A_j = H_{2j-1}+i H_{2j} \quad  \hbox{ with } \ H_{2j-1}, H_{2j} \in \SH
\quad  { for  } \ j = 1, \dots, m.$$
Then $\Lambda_k(\bA) \subseteq \IC^{m}$
can be identified with $\Lambda_k(H_1, \dots, H_{2m}) \subseteq \IR^{2m}$.
Thus,   we will focus on the joint rank $k$-numerical ranges
of self-adjoint operators in our discussion.

\section{Basic properties of $\Lambda_k(\bA)$}

\begin{proposition}   \label{2.1}
Suppose $\bA = (A_1, \dots, A_m) \in \SH^m$,
and $T = \(t_{ij}\)$ is an $m\times n$ real matrix.
If $B_j = \sum_{i=1}^m t_{ij} A_i$ for $j = 1, \dots, n$, then
$$  \{\ba T: \ba \in \Lambda_k(\bA)\}\subseteq \Lambda_k(\bB).$$
Equality holds if $\{A_1, \dots, A_m\}$ is linearly independent
and  $\span\{A_1, \dots, A_m\}= \span\{B_1, \dots, B_n\}$.
\end{proposition}

\it Proof. \rm The set inclusion follows readily from definitions.
Evidently, the equality holds if $n=m$ and $T$ is invertible.

Suppose $\{A_1, \dots, A_m\}$ is linearly independent and  span $\{A_1, \dots,
A_m\}=$ span $\{B_1, \dots, B_n\}$. First consider the special case when
$A_i=B_i$ for $1\le i\le m$. Then $T=\[I_m|T_1\]$ for some $m\times (n-m)$
matrix $T_1$. Let $\(b_1,\dots, b_n\)\in \Lambda_k(\bB)$. Then there exists a
rank $k$ orthogonal projection $P$ such that $PB_iP=b_iP$ for $i = 1, \dots,
n$. Therefore, we have $\(b_1,\dots,b_m\)\in  \Lambda_k(\bA)$ and for $1\le
j\le n$,
$$b_jP= PB_jP=P\(\sum_{i=1}^m t_{ij} A_i\)P=P\(\sum_{i=1}^m t_{ij}
B_i\)P=\(\sum_{i=1}^m t_{ij} b_i\)P\Ra b_j=\(\sum_{i=1}^m t_{ij} b_i\).$$
Therefore, $\(b_1,\dots, b_n\)=\(b_1,\dots,b_m\)T$.

For the general case, by applying a permutation, if necessary, we may assume
that $\{B_1,\dots, B_m\}$ is a basis of  $\span\{B_1, \dots, B_n\}$. Then
there exists an $m\times m$ invertible matrix $S=\(s_{i\,j}\)$
such that $A_j = \sum_{i=1}^m s_{ij} B_i$ for $j = 1, \dots, m$.
For $1\le j\le m$, we have
$$B_j=\sum_{i=1}^m t_{ij} A_i=\sum_{i=1}^m
t_{ij}\(\sum_{k=1}^m s_{ki}B_k\)=\sum_{k=1}^m \(\sum_{i=1}^m
s_{ki}t_{ij}\)B_k\,.$$
Therefore,  $\sum_{i=1}^m s_{ki}t_{ij}=\delta_{k\,j}$
and $ST=\[I_m|T_1\]$ for some $m\times (n-m)$ matrix $T_1$. Hence, we have
$$\Lambda_k(\bB)=\Lambda_k\(B_1,\dots,B_m\)\[I|T_1\]=
\Lambda_k\(B_1,\dots,B_m\)ST=\Lambda_k\(A_1,\dots,A_m\)T.$$
\vskip -.3in \qed

\medskip
In view of the above proposition, in the study of the geometric properties
of $\Lambda_k(\bA)$,
we may always assume that $A_1, \dots, A_m$ are linearly independent.

\begin{proposition} \label{2.2} Let $\bA = (A_1, \dots, A_m) \in \SH^m$,
and let $k < \dim \cH$.

{\rm (a)}
For any real vector ${\bmu} = (\mu_1, \dots, \mu_m)$,
$\Lambda_k(A_1 - \mu_1 I, \dots, A_m - \mu_m I) =
\Lambda_k(\bA) - {\bmu}$.

{\rm (b)} If $(a_1, \dots, a_m) \in \Lambda_k(\bA)$ then
$(a_1, \dots, a_{m-1}) \in \Lambda_{k}(A_1, \dots, A_{m-1})$.

{\rm (c)}
$\Lambda_{k+1}(\bA) \subseteq \Lambda_{k}(\bA)$.
\end{proposition}

\begin{remark} \label{2.3} \rm
By Proposition \ref{2.2} (a),
we can replace $A_j$ to $A_j - \mu_j I$ for $j = 1, \dots, m,$
without affecting  the geometric properties of
$\Lambda_k(A_1, \dots, A_m)$.
\end{remark}

Suppose $\dim \cH =n< 2k-1$ and
$A_1 = \diag(1,2,\dots, n)$. Then $\Lambda_k(A_1) = \emptyset$.
By Proposition \ref{2.2} (c), we see that $\Lambda_k(\bA) = \emptyset$
for any $A_2, \dots, A_m$. Thus, $\Lambda_k(\bA)$ can be empty
if $\dim \cH$ is small.
However,  a result of Knill, Laflamme  and Viola  \cite{KLV} shows that
$\Lambda_k(\bA)$ is non-empty if $\dim \cH$ is sufficiently large. By
modifying the proof of Theorem 3 in \cite{KLV}, we can get a slightly better
bound in the following proposition. The proof given here is essentially the
same as that of Theorem 3 and 4 in \cite{KLV}, except the choice of $\x_1$. We
include the details here for completeness.

\begin{proposition}  \label{2.4}
Let $\bA \in \SH^m$. For $m\ge 1$ and $ k> 1$.
If $\dim \cH =n\ge (k-1)(m+1)^2$,
then $\Lambda_k(\bA) \ne \emptyset$.
\end{proposition}

\it Proof. \rm
 We may assume that $\dim \cH = n = (k-1)(m+1)^2$.
Otherwise, replace each $A_j$ by $U^*A_jU$ for some $U$ such that $U^*U =
I_n$.
Let $q=(m+1)(k-1)+1$. Choose an eigenvector $\x_1$ of $A_1$ with $\|\x_1\|=1$.
Then choose a unit vector $\x_2$ orthogonal to $\x_1, A_2\x_1,\dots,A_m\x_1$.
By the assumption on $n$, we can choose an orthonormal set $\{\x_1,\x_2,\dots,
\x_q\} $ of $q$ vectors in $\IC^n$ such that for $1<r\le q$, $\x_r$ is
orthogonal to $A_j\x_i$ for all $1\le i<r$ and $1\le j\le m$. Let $X$ be the
$n\times q$ matrix with $\x_i$ as the $i$-th column. Then $X^*A_jX$ is a
diagonal matrix for $1\le j\le m$. By Tverberg's Theorem \cite{Tv}, we can
partition the set $\{i:1\le i\le q\}$ into $k$ disjoint subset $R_j$, $1\le
j\le k$ such that  $R =\cap_{j=1}^k\conv \{\la \bA\x_i,\x_i\ra:i\in R_j\}\ne
\emptyset$. Suppose $\ba\in R$. Then there exist non-negative numbers
$t_{i\,j}$, $1\le j\le k$, $i\in R_j$ such that for all $1\le j\le k$,
$\sum_{i\in R_j}t_{i\,j}=1$ and $\sum_{i\in R_j}t_{i\,j}\la
\bA\x_i,\x_i\ra=\ba$. Let $\y_j= \sum_{i\in R_j}\sqrt{t_{i\,j}}\x_i$ for
$1\le j\le k$. Then $\{\y_1,\dots,\y_k\}$ is orthonormal and $\la
\bA\y_j,\y_j\ra=\ba$ for all $1\le j\le k$.
\qed

\begin{proposition} \label{2.5}
Suppose $\bA \in \SH^m$ and $1 \le r < k \le \dim \cH$.
Let $\cV_r$ be the set of operator $X: \cH_1^\perp \rightarrow \cH$
such that $X^*X = I_{\cH_1^\perp}$ for an $r$-dimensional subspace $\cH_1$ of
$\cH$. Then
\begin{equation}
\label{inclu}
\Lambda_k(\bA) \subseteq
\cap \{ \Lambda_{k-r}(X^*A_1X, \dots, X^*A_mX): X \in \cV_r\}
\end{equation}
and
\begin{equation}\label{inclu2}
\conv \Lambda_k(\bA) \subseteq
\conv\left(\cap \{ \Lambda_{k-r}(X^*A_1X, \dots, X^*A_mX):
X \in \cV_r\}\right).
\end{equation}
\end{proposition}

\it Proof. \rm
Suppose $(a_1, \dots, a_m) \in \Lambda_k(\bA)$.
Let $\cH_2$ be a $k$-dimensional subspace and
$V: \cH_2 \rightarrow \cH$ such that $V^*V = I_{\cH_2}$
and  $V^*A_jV  = a_j I_k$ for $j = 1, \dots, m$.
Let $X\in \cV_r$ and $X^*X = I_{\cH_1^\perp}$ for an $r$-dimensional
subspace $\cH_1$ of $\cH$.
Then $\cH_0 =X^*\(V\(\cH_2\) \cap X\(\cH_1^\perp\)\)$
has dimension at least $s = k-r$. Let
$U: \cH_0 \hookrightarrow \cH$ be given by $Ux=x$ for all $x\in \cH_0$.
Then we have $U^*U = I_{\cH_0}$ and
$U^*\(X^*A_jX\)U  = a_j I_{\cH_0}$ for $j = 1, \dots, m$.
Thus, (\ref{inclu}) holds, and the inclusion
(\ref{inclu2}) follows.
\qed

Proposition \ref{2.5} extends \cite[Proposition 4.8]{LPS2}
corresponding to the case when $m = 2$. In such a case,
the set inclusion (\ref{inclu}) becomes a set equality if
$\dim \cH < \infty$ or if $(A_1,A_2)$ is a commuting pair, i.e.,
$A_1+iA_2$ is normal; see \cite[Corollary 4.9]{LPS2}.
The following example shows that the
set equality in  (\ref{inclu}) may not hold even in
the finite dimensional case if $m\ge 3$.

\begin{example} \label{2.6} \rm Let
$B_1 = \(\begin{array}{rr}1 & 0 \\ 0 & -1\end{array}\)$
$B_2 =  \(\begin{array}{rr}0 & 1 \\ 1 & 0\end{array}\)$
$B_3 =  \(\begin{array}{rr}0 & i \\ -i & 0\end{array}\)$.
For $k > 1$, let $A_j = B_j \otimes I_k$ for $j = 1, 2,3$.

(a) We have $\Lambda_k(A_1, \dots, A_m) =
\Lambda_1(B_1, B_2, B_3)=\{\ba\in\IR^3:\|\ba||=1\}$,
which  is not convex.

(b)
If $r=k-1$ and $X\in \cV_r$ and $X^*X=I_{\cH_1^\perp}$
for an $r$-dimensional subspace $\cH_1$, then
$\dim \cH_1^\perp =2k-r=k+1\ge 3$ so that
$$\Lambda_{k-r}\(X^*A_1X, X^*A_2X, X^*A_3X\)
=\Lambda_1\(X^*A_1X, X^*A_2X, X^*A_3X\)$$
is convex \cite{LP}.

\medskip\noindent
Consequently,
$\cap \{ \Lambda_{k-r}\(X^*A_1X, X^*A_2X, X^*A_3X\): X \in \cV_r\}$
is convex and cannot be equal to $\Lambda_k(A_1,A_2,A_3)$.

\medskip
For $m>3$, we can take $A_1$, $A_2$, $A_3$ as above and $A_j=0_{2k}$
for $3< j\le m$. Then we have
$$\Lambda_k(A_1,A_2,A_3) \ne
\cap \{ \Lambda_{k-r}(X^*A_1X, \dots, X^*A_mX): X \in \cV_r\}\,.$$
\end{example}

To verify (a), suppose
$U = \begin{pmatrix} U_1 \cr U_2 \cr\end{pmatrix}$
is such that $U_1, U_2 \in M_k$,
$U^*U  = U_1^*U_1 + U_2^*U_2 = I_k$ and $U^*A_jU = a_j I_k$.
Then $U_1^*U_1 - U_2^*U_2 = a_1 I_k$.  It follows that
$U_1^*U_1 = (1+a_1)I_k$ and $U_2^*U_2 = (1-a_1)I_k$
Thus, $U_1U_1^* = (1+a_1)I_k$ and $U_2U_2^* = (1-a_1)I_k$.
As a result,
$$U_i^*U_jU_j^*U_i = (1+a_1)(1-a_1)I = U_i^*U_iU_j^*U_j \quad
\hbox{ for } (i,j) \in \{(1,2), (2,1)\},$$
and
\begin{eqnarray*}
(a_1^2+a_2^2+a_3^2)I_k
&=& \sum_{j=1}^3 (U_j^* A_j U_j)^2    \\
&=& (U_1^*U_1-U_2^*U_2)^2 + (U_1^*U_2+U_2^*U_1)^2 +
(iU_1^*U_2-iU_2^*U_1)^2  \\
&=& (U_1^*U_1+U_2^*U_2)^2  \\
&=& I_k.
\end{eqnarray*}
Thus, $\Lambda_k(A_1,A_2,A_3) \subseteq \{(a_1, a_2, a_3)\in \IR^3:
a_1^2 + a_2^2+a_3^2 = 1\}.$

Conversely, suppose
$(a_1, a_2, a_3)\in \IR^3$ such that $ a_1^2 + a_2^2+a_3^2 = 1$.
Let
$$(\alpha,\beta) = \begin{cases}
 (0,1) &  if \  a_1 = -1,   \\    & \\
\dfrac{\left(1+a_1, a_2-ia_3\right)}
{\sqrt {2(1+a_1)}}&  otherwise.
\end{cases}
$$
Let $U = \begin{pmatrix} \alpha I_k \cr \beta I_k \cr\end{pmatrix}$.
Direct computation shows that $U^*A_jU=a_jI_k$ for $1\le j\le 3$.

By a similar argument or putting $k = 1$,
we see that $\Lambda_1(B_1,B_2,B_3)$ has the same form.
\qed

\medskip
It is natural to ask if the set equality in (\ref{inclu2}) can hold for
for $m > 2$. Also,
$$\conv\left(\cap \{ \Lambda_{k-r}(X^*A_1X, \dots, X^*A_mX):
X \in \cV_r\}\right) $$
$$\subseteq
\cap \{\conv  \Lambda_{k-r}(X^*A_1X, \dots, X^*A_mX):
X \in \cV_r\}.$$
It is interesting to determine whether the two sets are equal.

\section{Geometric properties of $\Lambda_k(\bA)$}
\setcounter{equation}{0}

Let $\bA = (A_1, \dots, A_m) \in \SH^m$.
It is known that $\Lambda_k(\bA)$ is always convex
if $m \le 3$ unless $(\dim \cH, m) = (2,3)$.
If $(\dim \cH, m) = (2,3)$ or $n>1,\ m \ge 4$,
there are examples $\bA \in \SH^m$ such that $\Lambda_1(\bA)$
is not convex.  Furthermore, for any $A_1, A_2, A_3\in \SH$
such that $\span\{I, A_1, A_2, A_3\}$ has dimension 4, there
is always an $A_4 \in \SH$ for which
$\Lambda_1(A_1,\dots, A_4)$ is not convex; see \cite{LP}.
Nevertheless,  for $m \ge 4$
we will show that $\Lambda_k(\bA)$ is always star-shaped
if $\dim \cH$ is sufficiently large.
Moreover, it always contains the convex hull of
$\Lambda_{\hat k}(\bA)$ for  $\hat k = (m+2)k$. If $k = 1$,
the result can be
further improved.
We begin with the following result.

\begin{theorem} \label{3.1}
Let $\bA = (A_1, \dots, A_m) \in \SH^m$ and $k$ be a positive integer.
If $\Lambda_{\hat k}(\bA) \ne \emptyset$
for some $\hat k \ge (m+2)k$, then
$\Lambda_k(\bA)$ is star-shaped and contains the convex
subset $\conv\Lambda_{\hat k}(\bA)$ so that every element
in $\conv\Lambda_{\hat k}(\bA)$ is
a star center of $\Lambda_k(\bA)$.
\end{theorem}

Note that $\Lambda_{\hat k}(\bA)$ may be empty
if $\dim \cH$ is small relative to $\hat k$.
Even if $\Lambda_{\hat k}(\bA)$  is non-empty, it may be much smaller
than its convex hull; for example,
see Example \ref{2.6}. So, the conclusion in Theorem \ref{3.1}
is rather remarkable.

\it Proof. \rm We may assume $\ba = \b0 \in \Lambda_{\hat k}(\bA)$. Then
there exists $Y$ such that $Y^*Y = I_{(m+2)k}$ and
$Y^*A_jY =\b0$ for all $\le j\le m$.

Let $\bb\in \Lambda_k(\bA)$. Then there exists $X$   such that
$X^*X = I_k$ and
$$(X^*A_1X, \dots, X^*A_mX) = (b_1 I_k, \dots, b_m I_k).$$
Suppose $\cX$ and $\cY$ are the range spaces of $X$ and $Y$,
respectively. Then we have
$$\dim \(\cY\cap\(\cX+A_1(\cX)+\cdots+A_m(\cX)\)^{\perp}\)
\ge \dim \cY-\dim \(\cX+A_1(\cX)+\cdots+A_m(\cX)\)\ge k.$$
Let $\cY_{1}$ be an $k$-dimensional subspace of
$\cY\cap\(\cX+A_1(\cX)+\cdots+A_m(\cX)\)^{\perp}$ and
$\cY_{2}=\cY\cap \(\cX+\cY_1\)^{\perp}$.  Set $Z=\[X|Y_1|Y_{2}\]$,
where $Y_i$ has columns forming an orthonormal basis of $\cY_i$ for $i=1,2$.
Then we have  $Z^*Z = I_{(m+2)k}$ and for $1\le j\le m$,
$Z^*A_jZ$ has the form
$$\(  \begin{array}{ccc} b_j I_k & \b0_k&*  \\
\b0_k&\b0_k&*\\
* & *&* \end{array}\).$$
Let
$C_j = b_j I_k \oplus \b0_k$. For $t \in [0,1]$, we have
$$t(b_1, \dots, b_m) \in \Lambda_k(C_1, \dots, C_m)
\subseteq \Lambda_k(Z^*A_1Z,\dots, Z^*A_mZ) \subseteq \Lambda_k(\bA).$$

Clearly, $\Lambda_{\hat k}(\bA)\subseteq \Lambda_k(\bA)$.
Since every element of $\Lambda_{\hat k}(\bA)$
is a star center of $\Lambda_k(\bA)$ and the set of star centers of a star-shaped set is convex, we see that
every element in $\conv\Lambda_{\hat k}(\bA)$ is a star center of
$\Lambda_k(\bA)$. Hence, $\conv \Lambda_{\hat k}(\bA)\subseteq \Lambda_k(\bA)$.
\qed

If $\dim\cH$ is finite, then $\Lambda_k(\bA)$ is always closed.
But this may not be the case if $\dim\cH$ is infinite.
Using Theorem \ref{3.1}, we can prove the star-shapedness
of $\cl(\Lambda_k(\bA))$.

\begin{corollary} \label{3.2}
Let $\bA = (A_1, \dots, A_m) \in \SH^m$ and $k$ be a positive integer.
If $\cl(\Lambda_{\hat k}(\bA)) \ne \emptyset$
for some $\hat k \ge (m+2)k$, then
$\cl(\Lambda_k(\bA))$ is star-shaped and contains the convex
subset $\conv\cl(\Lambda_{\hat k}(\bA))$ so that every element
in $\conv\cl(\Lambda_{\hat k}(\bA))$ is
a star center of $\cl(\Lambda_k(\bA))$.
\end{corollary}

\it Proof. \rm Suppose $\ba \in  \cl(\Lambda_{\hat k}(\bA))$
and $\bb \in \Lambda_k(\bA)$. Then for every $\varepsilon$ there
is $\tilde \ba = (\tilde a_1, \dots, \tilde a_m) \in \Lambda_{\hat k}(\bA)$
such that $\ell_1(\tilde \ba - \ba) < \varepsilon$.
By Theorem \ref{3.1}, we see that the line segment joining
$\tilde \ba$ and $\bb$ lies in $\Lambda_k(\bA)$.
Consequently, the line segment joining $\ba$ and $\bb$ lies in
$\cl(\Lambda_k(\bA))$. The proof of the last assertion is similar to that
of Theorem \ref{3.1}.
\qed

It is easy to see that a star center of $\Lambda_k(\bA)$
is also a star center of  $\cl(\Lambda_k(\bA))$. However,
the converse may not hold.
The following example from \cite[Example 3.2]{LP2}
illustrates this.

\medskip
\begin{example}\label{3.3}
Consider $\cH = \ell^2$ with canonical basis $\{e_n: n \ge 1\}$.
Let $\bA = (A_1, \dots, A_4)$ with
$A_1=\diag(1,0,1/3,1/4,....), A_2 =\diag(1,0)\oplus 0,$
$$A_3=\begin{pmatrix} 0&1\cr 1&0\end{pmatrix} \oplus 0 \quad \hbox{ and }
\quad A_4=\begin{pmatrix} 0&i\cr -i&0\end{pmatrix} \oplus 0.$$
Then $(1,1,0,0)\in W(\bA)$ and
$(0,0,0,0) \in  W(\bA)\cap W_e(\bA)$
is a star-center of $\cl(W(\bA))$. However,
$(1/2,1/2,0,0) \notin W(\bA)$ so that $(0,0,0,0)$
is not a star-center of $W(\bA)$. In fact,
$W(\bA)$ is not convex even though $\cl(W(\bA))$ is convex.
\end{example}

By Proposition \ref{2.4}, we see that $\Lambda_{\hat k}(\bA)$ is
non-empty if $\dim\cH$ is sufficiently large. So, $\Lambda_k(\bA)$
is star-shaped and contains a convex set.
The same comment also holds for $\cl(\Lambda_k(\bA))$.
More specifically, we have the following.

\begin{theorem} \label{3.4}
Let $\bA = (A_1, \dots, A_m) \in \SH^m$. If
$\dim \cH \ge ((m+2)k-1)(m+1)^2$, then both $\Lambda_k(\bA) $ and
$\cl(\Lambda_k(\bA))$ are star-shaped.
\end{theorem}

In Theorem \ref{3.7}, we will show that
the classical joint numerical range is star-shaped
with a much milder restriction on $\dim\cH$ comparing with
that in Theorem \ref{3.4}.
To demonstrate this, we need two related results.

\begin{proposition} \label{3.5} Suppose $\dim\cH = n$ and
$\bA = (A_1, \dots, A_m) \in \SH^m$ is such that
$\{A_1, \dots, A_m\}$ is linearly independent.
Assume that
$\b0 \in \Lambda_{n-1}(\bA)$, i.e., there is a basis
such that $A_j$ has operator matrix
$\begin{pmatrix}* & *  \cr * & 0_{n-1}\end{pmatrix}$
for $j = 1, \dots, m$.

{\rm (a)} If $m = 2n-1$, then there is an invertible $S \in M_m(\IR)$
such that
$$\Lambda_1(\bA) =
\left\{(1+u_1,u_2,\dots, u_m)S:
u_1, \dots, u_m \in \IR, \sum_{j=1}^m u_j^2 = 1\right\}$$
so that $\Lambda_1(\bA)$ is not star-shaped.

{\rm (b)} If $m < 2n-1$, then $\Lambda_1(\bA)$ is star-shaped with
$\b0$ as a star center.
\end{proposition}

\it Proof. \rm (a)
If $m = 2n-1$, there is
an invertible $m\times m$ real matrix $T = (t_{ij})$
such that for $B_j = \sum_{i=1} t_{ij} A_i$ for $j = 1, \dots, m$
with
$$\bB =
(B_1, \dots,  B_m) = (E_{11}, E_{12}+E_{21}, -iE_{12}+iE_{21},
\dots,  -iE_{1n}+iE_{n1}).$$
Note that every unit vector $x \in \IC^n$ has the form
$\x = \mu (\cos t, \sin t(v_2+iv_3), \dots, \sin t(v_{m-1}+iv_m))^t$
such that $|\mu|  = 1$, $t \in [0,\pi/2]$,
and $v_2, \dots, v_m \in \IR$ with $\sum_{j=2}^m v_j^2 = 1$.
We have
\begin{eqnarray*}
\langle \bB \x, \x \rangle
&=& ((1+\cos(2t))/2, v_2\sin (2t), v_2 \sin(2t), \dots, v_m \sin(2t)) \\
&=& (1+u_1,u_2, \dots, u_m)D,
\end{eqnarray*}
where  $D = [1/2] \oplus I_{n-1}$,
$u_1 = \cos(2t)$ and $u_j = v_j \sin(2t)$ for $j =2, \dots, m$.
It follows that
$$\Lambda_1(\bB) =
\left\{(1+u_1, u_2, \dots, u_m)D:
u_1, \dots, u_m \in \IR, \sum_{j=1}^m u_j^2 = 1\right\}.$$
By Proposition \ref{2.1},
$\Lambda_1(\bA) = \{ \bb T^{-1}: \bb \in \Lambda_1(\bB)\}$.
The result follows.

(b)
Now, suppose $m < 2n-1$. By adding more $A_j$, if necessary,
we only need to consider the case when $m = 2n-2$.
Let $\v_j$ be the row vector obtained by
removing the first entry of the first row of $A_j$
for $j = 1, \dots, m$.

\medskip\noindent
{\bf Case 1}
Suppose $\span\{\v_1, \dots, \v_m\}$ has real dimension $m-1 = 2n-3$.
Then there is a unitary matrix of the form $U = [1] \oplus U_0 \in M_n$
such that the $(1,n)$ entry of $U^*A_jU$ is real  for $j = 1, \dots, m$.
Hence, there is
an invertible $m\times m$ real matrix $T = (t_{ij})$
such that for $B_j = \sum_{i=1} t_{ij} A_i$ for $j = 1, \dots, m$
with
\begin{eqnarray*}
\bB &=&  (B_1, \dots,  B_m)\\
&=& (E_{11}, E_{12}+E_{21}, -iE_{12}+iE_{21},
\dots,  -iE_{1,n-1}+iE_{n-1,1}, E_{1n}+E_{n1}).
\end{eqnarray*}
Suppose $\bb\in\Lambda_1(\bB)$, $\bb\ne\b0$.
Then there exists a unit vector
$\x = \mu (u_0, u_1+iu_2, \dots, u_{m-1}+iu_m)^t$
such that  $|\mu| = 1$ and $u_0 > 0$, with
$$\bb=\langle \bB \x, \x \rangle =
u_0(u_0, 2u_1, \dots, 2u_{m-1}).$$
For any $t \in (0,1)$, we can choose a unit vector of the form
$\x_t = \sqrt{t} (u_0, u_1 + iu_2,  \dots,  u_{m-1}+i\tilde u_m)^t$
with $t\tilde u_m^2 = 1-\sum_{j=0}^{m-1} t u_j^2$
so that
$$\langle \bB \x_t, \x_t \rangle =
t\langle \bB \x, \x \rangle.$$

\medskip\noindent
{\bf Case 2.} Suppose $\{\v_1, \dots, \v_m\}$
has real dimension $m = 2n-2$.
Then there is an invertible $m\times m$ real matrix $T = (t_{ij})$
such that for $B_j = \sum_{i=1} t_{ij} A_i$ for $j = 1, \dots, m$
with
\begin{eqnarray*}
\bB &=&  (B_1, \dots,  B_m) \\
& = & (a_1 E_{11}, \dots, a_m E_{11}) +
(E_{12}+E_{21}, iE_{12}+iE_{21},
\dots,  E_{1n}+E_{n1}, -iE_{1n}+iE_{n1})
\end{eqnarray*}
with $a_1, \dots, a_m \in \IR$.
Suppose $\bb\in\Lambda_1(\bB)$, $\bb\ne\b0$.
Then there exists a unit vector
$\x = \mu (u_0, u_1+iu_2, \dots, u_{m-1}+iu_m)^t$
such that  $|\mu| = 1$ and $u_0 > 0$, with
$$\bb=\langle \bB \x, \x \rangle =
u_0(a_1 u_0 + 2u_1, a_2u_0 + 2u_2, \dots, a_m u_0 + 2u_{m}).$$
For any $t \in (0,1)$, consider a vector of the form
$$\x_\xi = (\xi u_0, w_1+
i w_2, w_3+iw_4, \dots, w_{m-1}+iw_m)^t,$$
where $\xi \ge t$ and $w_j = a_ju_0(t-\xi^2)/(2\xi) + tu_j/\xi$
for $j = 1, \dots, m$.
Then
$$\xi u_0 (a_j \xi u_0 + 2w_j)
= tu_0(a_j u_0 + 2u_j), \qquad j = 1, \dots, m,$$
so that
$$\langle \bB \x_\xi, \x_\xi \rangle =
t \langle \bB \x, \x \rangle.$$
If $\xi = \sqrt{t}$, then
$\x_\xi = \sqrt t (u_0, u_1+iu_2, \dots, u_{m-1}+iu_m)^t$
has norm less than 1;
if $\xi \rightarrow \infty$, then
$\|\x_\xi\| \ge |\xi u_0| \rightarrow \infty$.
Thus, there is $\xi > t$ such that   $\x_\xi$ is a unit vector
satisfying
$\langle \bB \x_\xi, \x_\xi \rangle =
t \langle \bB \x, \x \rangle.$
So,  $\Lambda_1(\bB)$ is star-shaped with $\b0$ as a star center.
By Proposition \ref{2.1},
$\Lambda_1(\bA) = \{ \bb T^{-1}: \bb \in \Lambda_1(\bB)\}$.
The result follows.
\qed

\begin{theorem} \label{3.6} Let $\bA = (A_1,\dots, A_m) \in \SH^m$.
If $\Lambda_{\hat k}(\bA) \ne \emptyset$
for some $\hat k > (m+1)/2$,
then $\Lambda_1(\bA)$ is star-shaped and contains
$\conv \Lambda_{\hat k}(\bA)$ such that every element in
$\conv \Lambda_{\hat k}(\bA)$ is a star center of
$\Lambda_1(\bA)$.
\end{theorem}

\it Proof. \rm We may assume $\ba = \b0 \in \Lambda_{\hat k}(\bA)$
with $\hat k > (m+1)/2$.
Suppose $\x \in \cH$ is a unit vector and
$\bb  = \la \bA \x, \x \ra \in \Lambda_1(\bA)$.
Suppose $X$ is such that $X^*X = I_{\hat k}$ and
$X^* A_j X = 0_{\hat k}$ for $j = 1, \dots, m$.
Let $Y$ be such that $Y^*Y = I_{\hat k+1}$, and
the range space of $Y$ contains the range space of $X$ and $\x$.
Suppose $\bB = (B_1, \dots, B_m) = (Y^*A_1Y, \dots, Y^*A_mY)$.
Then we may assume that $B_j$ has the form
$\begin{pmatrix}* & *  \cr * & 0_{\hat k}\end{pmatrix}$ for $j = 1, \dots, m$.
Clearly, $\span\{B_1, \dots, B_m\}$ has dimension at most $m < 2\hat k-1$.
By Proposition \ref{3.5}, the line segment joining $\b0$ and
$\bb$ lies entirely in $\Lambda_1(\bB) \subseteq \Lambda_1(\bA)$.
Thus, $\b0$ is a star center of $\Lambda_1(\bA)$. Since the
set of star centers of $\Lambda_1(\bA)$ is convex, we see that
every element in $\conv\Lambda_{\hat k}(\bA)$ is a star center
of $\Lambda_1(\bA)$. \qed

\begin{theorem} \label{3.7} Let $\bA = (A_1,\dots, A_m) \in \SH^m$.
If  $\dim\cH \ge \[\dfrac{m+1}2\] (m+1)^2$,
then $\Lambda_1(\bA)$ is star-shaped.
\end{theorem}

\it Proof. \rm Let $\hat k= \[\dfrac{m+1}2\] +1> \dfrac{m+1}2$.
Then $\dim \cH \ge (\hat k - 1)(m+1)^2$ and $\Lambda_{\hat k}(A) \ne
\emptyset$ by Proposition \ref{2.4}.
The result then follows from Theorem \ref{3.6}. \qed

\section{Results on $\Lambda_\infty(\bA)$}

In this section, we always assume that $\cH$ has infinite dimension.
Denote by
$\cP_\infty$ the set of infinite rank orthogonal projections
in $\SH$, and let
$$\Lambda_\infty(\bA)
= \{ (\gamma_1, \dots, \gamma_m) \in \IR^m: \hbox{ there is } P \in \cP_\infty
\hbox{ such that } PA_iP = \gamma_i P \ \mbox{ for all } 1\le i\le m\}\,$$
for $\bA \in \SH^m$.

By the result in Section 3, we have the following.

\begin{proposition}\label{4.1}
Suppose $\bA \in \SH^m$, where $\cH$ is infinite-dimensional.
Then $\Lambda_k(\bA)$ is star-shaped for each positive integer $k$.
Moreover, if $\ba \in \Lambda_\infty(\bA)$,
then $\ba$ is a star center for $\Lambda_k(\bA)$ for every positive
integer $k$.
\end{proposition}

When $m = 2$, it was conjectured in \cite{R} and confirmed in
\cite{LPS2} that
$$\Lambda_\infty(A_1,A_2) = \bigcap_{k \ge 1} \Lambda_k(A_1,A_2);$$
in \cite[Theorem 4]{AS}, it was proven that
$$\Lambda_\infty(A_1,A_2)
=\cap\{W(A_1+F_1,A_2+F_2): F_1, F_2 \in \SH \cap \FH\}.$$
In the following, we extend the above results
to $\Lambda_\infty(A_1, \dots, A_m)$ for $m > 2$. Moreover, we show that
$\Lambda_\infty(\bA) =
\bigcap_{k \ge 1} S_k(\bA)$,
where $S_k(\bA)$ is the set of star centers of $\Lambda_k(\bA)$.
Hence, $\Lambda_\infty(\bA)$ is always convex.

\begin{theorem}\label{4.2}
Suppose $\bA \in \SH^m$, where $\cH$ is infinite-dimensional.
For each $k\ge 1$, let $S_k(\bA)$ be the set of star-center of
$\Lambda_k(\bA)$. Then
\begin{equation}\label{eqn4}\Lambda_\infty(\bA)=
 \cap_k S_k(\bA)= \cap_k\Lambda_k(\bA)
= \cap\{ W(\bA+\bF): \bF \in\SH^m \cap \FH^m  \}\,.\end{equation}
Consequently, $\Lambda_\infty(\bA)$ is  convex.
\end{theorem}

\it Proof. \rm It follows from definitions and Theorem \ref{3.1} that
$$\Lambda_\infty(\bA)\subseteq
 \cap_k S_k(\bA)\subseteq \cap_k\Lambda_k(\bA)\,.$$
We are going to prove that
$\cap_k\Lambda_k(\bA)\subseteq \cap\{ W(\bA+\bF):
\bF \in\SH^m \cap \FH^m  \}$.
Suppose $\bF=(F_1,\dots,F_m) \in\SH^m \cap \FH^m$
and  $K=\sum_{i=1}^m\rank\(F_i\)+1$.
Let $\bmu=(\mu_1,\dots,\mu_m)\in  \cap_k\Lambda_k(\bA)$. Then there exists
a rank $K$ orthogonal projection $P$ such that $PA_jP=\mu_jP$
for $1\le j\le m$. Let
$$\cH_0=\range P\cap \ker F_1\cap \ker F_2\cap\cdots\cap\ker F_m
= \range P\cap \(\range F_1+\range F_2+\cdots+\range F_m\)^\perp\,.$$
Then dim $\cH_0\ge 1$. Let $\x$ be  a unit vector  in $\cH_0$. Then  we have
$\la \(\bA+\bF\)\x,\x\ra=\la \bA\x,\x\ra=\bmu$.  Therefore,
$\bmu\in  W(\bA+\bF)$. Hence, we have
$\cap_k\Lambda_k(\bA)\subseteq \cap\{ W(\bA+\bF) :
\bF \in\SH^m \cap \FH^m  \}$.

Next, we prove that $\cap\{ W(\bA+\bF) :
\bF \in\SH^m \cap \FH^m  \}\subseteq \Lambda_\infty(\bA)$.
Suppose
$$\bmu\in  \cap\{ W(\bA+\bF) : \bF \in\SH^m \cap \FH^m  \}.$$
By Remark \ref{2.3}, we may assume that $\bmu=\b0$.
Then $\b0\in W(\bA)$ and there exists a unit vector $\x_1$
such that $\la \bA \x_1,\x_1\ra=\b0$. Suppose we have chosen an orthonormal
set of vectors $\{\x_1,\dots,\x_n\}$ such that $\la \bA \x_i,\x_j\ra =\b0$
for all $1\le i,\ j\le n$. Let $\cH_0$ be the subspace spanned by
$$\{\x_i:1\le i\le n\}\cup\{A_j\x_i:1\le i\le n,\ 1\le j\le m\}
\,$$
and $P$ the orthogonal projection of $\cH$ onto $\cH_0$.
Suppose
$$\bB=\(\left.(I-P)A_1(I-P)\right|_{\cH_0^\perp},\dots,
\left.(I-P)A_m(I-P)\right|_{\cH_0^\perp}\).$$
Let $\bb=\(b_1,\dots,b_m\) $ be a star-center of $W(\bB)$. Then
$$\bb I_{\cH_0}\oplus\bB=\(b_1P+(I-P)A_1(I-P),\dots,b_mP+(I-P)A_m(I-P)\)= \bA+\bF$$
for some $\bF\in\SH^m \cap \FH^m $. Therefore,
$\b0\in   W(\bb I_{\cH_0}\oplus\bB) $.
 Hence, there exists  a unit vector $\x\in \cH$ such that
$\b0=\la (\bA +\bF)\x,\x\ra$. Let $\x=\y+\z$,
where $\y\in \cH_0$ and $\z\in \cH_0^\perp$.
Then $\|\y\|^2+\|\z\|^2=\|\x\|^2=1$.
If  $\z=\b0$, then $\b0=\bb \in W(\bB) $. If
$\z\ne \b0$, then by Proposition \ref{4.1}, we have
$$\b0=\la (\bA  +\bF)\x,\x\ra
=\|\y\|^2\bb+\|\z\|^2\la \bB\(\dfrac{\z}{\|\z\|}\),
\(\dfrac{\z}{\|\z\|}\)\ra\in W(\bB)$$
So there exists a unit vector $\x_{n+1}\in \cH_0^\perp$ such that
$$\b0=\la (\bA +\bF)\x_{n+1},\x_{n+1}\ra=\la \bB \x_{n+1},\x_{n+1}\ra=\la \bA
\x_{n+1},\x_{n+1}\ra$$
Hence, inductively,  we can choose an orthonormal sequence of vectors
$\{\x_n\}_{n=1}^\infty$  such that
$$\la \bA \x_i,\x_j\ra=\b0\ \mbox{ for all }i,\ j\,.$$
Thus, we have
$$\Lambda_\infty(\bA)\subseteq \cap_k S_k(\bA) \subseteq
\cap_k\Lambda_k(\bA)\subseteq \cap\{ W(\bA+\bF) :
\bF \in\SH^m \cap \FH^m  \}\subseteq \Lambda_\infty(\bA).$$
Since $S_k(\bA)$ is convex for all $k\ge 1$, the last statement follows.
\qed

The last equality in (\ref{eqn4}) establishes a relationship between
$\Lambda_\infty(\bA)$ and the joint numerical ranges of finite rank
perturbation of $\bA$. The following result gives an extension.

\begin{theorem}\label{4.3}
Suppose $\bA = (A_1,\dots, A_m) \in \SH^m$.
Let  $n\ge 1$ and $\bF_0\in \SH^m\cap\FH^m$.
Then the following sets are equal.
\begin{itemize}
\item[a)] $ \Lambda_\infty(\bA)$.
\item[b)] $ \Lambda_\infty\(\bA+\bF_0\)$.
\item[c)] $  \cap\{\Lambda_n(\bA+\bF):\bF\in \FH^m\cap\SH^m\}$.
\item[d)] $\cap_{k\ge 1}\(\cap\{\Lambda_k(\bA+\bF):\bF\in \FH^m\cap\SH^m\}\)$.
\end{itemize}
\end{theorem}

\it Proof. \rm By Theorem \ref{4.2}, we have
$$ \begin{array}{rl}& \Lambda_\infty\(\bA+\bF_0\)\\
=& \cap\{ W(\bA+\bF_0+\bF): \bF \in\SH^m \cap \FH^m  \}\\
=& \cap\{ W(\bA+\bF): \bF \in\SH^m \cap \FH^m  \}\\
=&\Lambda_\infty\(\bA\).
\end{array}$$
This proves the equality of the sets in a) and b).
For the equality of the sets of a) and c),
let $\bF\in\FH^m\cap\SH^m$. Then we have
$$\Lambda_\infty\(\bA\)=\Lambda_\infty\(\bA+\bF\)
\subseteq \Lambda_n\(\bA+\bF\)\subseteq W(\bA+\bF).$$
It follows that
\begin{eqnarray*}
\Lambda_\infty\(\bA\)
&\subseteq& \cap\{\Lambda_n\(\bA+\bF\): \bF \in\SH^m \cap \FH^m  \} \\
&\subseteq& \cap\{W(\bA+\bF): \bF \in\SH^m \cap \FH^m  \}
=\Lambda_\infty\(\bA\).
\end{eqnarray*}

The equivalence of a) and d) follows immediately.\qed

Recall that $\cV_r$ is the set of $X:\cH_1^\perp\to \cH$ such that
$\dim \cH_1=r$ and $X^*X=I_{\cH_1^\perp} $.
Then $X^*\bA X$ is a compression of $\bA$ to $\cH_1^\perp$.
The next result is an analog to Theorem \ref{4.3} for $\Lambda_k\(X^*\bA
X\)$.

\begin{theorem}  \label{4.4}
Suppose $\bA = (A_1,\dots, A_m) \in \SH^m$.
Let  $n,\ r_0\ge 1$ and $X_0\in \cV_{r_0}$.
Then for $X^*\bA X = (X^*A_1X, \dots, X^*A_mX)$, the following sets are equal.
\begin{itemize}
\item[a)] $ \Lambda_\infty(\bA)$.
\item[b)] $ \Lambda_\infty\(X_0^*\bA X_0\)$.
\item[c)] $  \cap \{\Lambda_{n}(X^*\bA X):X\in \cup_{r\ge 1}\cV_r\}$.
\item[d)] $\cap_{k\ge 1}\( \cap
\{\Lambda_k(X^*\bA X):X\in \cup_{r\ge 1}\cV_r\}\)$.
\end{itemize}
\end{theorem}

\it Proof. \rm By (\ref{eqn4}) and (\ref{inclu}), we have
$$ \Lambda_\infty\(\bA\)
= \cap_k\Lambda_k(\bA)=\cap_k\Lambda_{k+r}(\bA)
\subseteq \cap_k\Lambda_{k}\(X_0^*\bA X_0\)
=\Lambda_\infty\(X_0^*\bA X_0\)\subseteq\Lambda_\infty(\bA).
$$
This proves the equality of the sets in a) and b). For the
equality of the sets of a) and c), we will first show that
\begin{equation}\label{eq6}
\cap\{\Lambda_1(X^*\bA X):X\in \cup_{r\ge 1}\cV_r\}\subseteq
\Lambda_\infty(\bA).
\end{equation}
Let $\bmu\in \cap\{\Lambda_1(X^*\bA X):X\in \cup_{r\ge 1}\cV_r\}$.
By Remark \ref{2.3}, we may assume that $\bmu=\b0$. Then
there exists a unit vector $\x_1$
such that $\la \bA \x_1,\x_1\ra=\b0$. Suppose we have chosen an orthonormal
set of vectors $\{\x_1,\dots,\x_N\}$ such that $\la \bA \x_i,\x_j\ra =\b0$
for all $1\le i,\ j\le N$. Let $\cH_1$ be the subspace spanned by
$$\{\x_i:1\le i\le N\}\cup\{A_j\x_i:1\le i\le N,\ 1\le j\le m\}\,$$
and $X:\cH_1^\perp\to \cH$ be given by $X(\v)=\v$ for all $\v\in\cH_1^\perp$.
Then $\b0\in \Lambda_1(X^*\bA X)$. So there exists a unit vector
$\x_{N+1}\in \cH_1^\perp$ such that
$$\b0=\la (X^*\bA X) \x_{N+1},\x_{N+1}\ra=\la\bA \x_{N+1},\x_{N+1}\ra\,.$$
Inductively, we can find an orthonormal sequence
$\{\x_i\}$ in $\cH$ such that $\la \bA\x_i,\x_j\ra=\b0$ for all $i,\ j$.
Hence, $\b0\in\Lambda_\infty(\bA)$.

To continue the proof of the equality of the sets of a) and c). Let
$X\in \cup_{r\ge 1}\cV_r$. Then we have
$$\Lambda_\infty\(\bA\) =\Lambda_\infty\(X^*\bA X\)\subseteq
\Lambda_n\(X^*\bA X\)\subseteq \Lambda_1\(X^*\bA X\).$$
It follows that
$$\Lambda_\infty\(\bA\)
\subseteq \cap\{\Lambda_\infty\(X^*\bA X\):X\in \cup_{r\ge 1}\cV_r \}
\subseteq \cap\{\Lambda_1(X^*\bA X):X\in \cup_{r\ge 1}\cV_r\}
\subseteq \Lambda_\infty\(\bA\).
$$
The equality of the sets  of a) and d) follows immediately.\qed

Recall that
the joint essential numerical range of $\bA\in \SH^m$ is defined by
$$W_e(\bA) = \cap \{\cl(W(\bA+\bF)): \bF \in \SH^m \cap \FH^m\}.$$
Using the last two theorems, we have the following.

\begin{corollary} \label{4.5}Let $\bA \in \SH^m$, where $\cH$ is
infinite dimensional.
Denote by $\tilde S_k(\bA)$ the set of star center of $\cl(\Lambda_k(\bA))$.
Then
$$W_e(\bA) = \bigcap_{k \ge 1} \cl\(\Lambda_k(\bA)\) =
\bigcap_{k\ge 1} \tilde S_k(\bA).$$
In addition, let  $n\ge 1$ and $\bF_0\in \SH^m\cap\FH^m$.
Moreover, let
$\bV$ be a finite dimensional subspace of $\cH$ and
$X_0: \bV^{\perp} \rightarrow \cH$ such that $X_0^*X_0 = I_{\bV^\perp}$.
Then   the following sets are equal.
\begin{itemize}
\item[a)] $ W_e(\bA)$.
\item[b)] $  W_e\(\bA+\bF_0\)$.
\item[c)] $  W_e\(X_0^*\bA X_0\)$.
\item[d)] $  \cap\{\cl\(\Lambda_n(\bA+\bF)\):\bF\in \FH^m\cap\SH^m\}$.
\item[e)] $  \cap\{\cl\(\Lambda_n(X^*\bA X)\):X\in \cup_{r\ge 1}\cV_r\}$.
\item[f)] $ \cap_{k\ge 1}\( \cap\{\cl\(\Lambda_k(\bA+\bF)\):
\bF\in \FH^m\cap\SH^m\}\)$.
\item[g)] $ \cap_{k\ge 1}\( \cap\{\cl\(\(\Lambda_k(X^*\bA X)\)\):
X\in \cup_{r\ge 1}\cV_r\}\)$.
\end{itemize}
\end{corollary}

Let $A \in \SH$ and $k$ be a positive integer. Denote by
$\cU_k$ the set of $X: \IC^k\rightarrow \cH$ such that
$X^*X = I_k$ and
$$\lambda_k(A) = \sup \{\min\sigma(X^*AX): X \in \cU_k \},$$
where $\sigma(B)$ is the spectrum of the operator $B$.
For $\bc=(c_1,\dots,c_m)\in\IR^m$ and
$\bA \in \SH^m$, let $\bc\cdot \bA=\sum_{i=1}^mc_iA_i$.
Define
$$\Omega_\infty(\bA) = \bigcap_{\bc \in \IR^m}
\left\{\ba \in \IR^m: \bc\cdot\ba
\le \lambda_k(\bc\cdot \bA) \mbox{ for all }k\ge 1\right\}.$$
The following extends \cite[Theorem 2.1]{LPS2}.

\begin{theorem} Let $\bA \in \SH^m$, where $\cH$ is infinite dimensional.
Then
$$\Lambda_\infty(\bA)\subseteq \Omega_\infty(\bA) =W_e(\bA). $$
\end{theorem}

\it Proof. \rm For $k\ge 1$, let
$$\Omega_k(\bA) = \bigcap_{\bc \in \IR^m}
\left\{\ba \in \IR^m: \bc\cdot\ba
\le \lambda_k(\bc\cdot \bA) \right\}.$$
Clearly, $\Omega_\infty(\bA)
=  \bigcap_{k\ge 1}\Omega_k(\bA)$. Suppose $k\ge 1$
and $\ba=(a_1,\dots,a_m)\in \Lambda_k(\bA)$.
Then there exists $P\in \cP_k$ such that $PA_jP=a_jP$ for all $1\le j\le m$.
For every $\bc \in \IR^m$, we have $P(\bc\cdot \bA)P=(\bc\cdot \ba)P$.
Hence, $$\bc\cdot \ba=\lambda_k(\bc\cdot P\bA P)=\lambda_k(P(\bc\cdot \bA) P)\le
\lambda_k(\bc\cdot  \bA).$$
Therefore, $\Lambda_k(\bA)\subseteq \Omega_k(\bA)$. Since $\Omega_k(\bA)$
is closed, we have $\cl(\Lambda_k(\bA))\subseteq \Omega_k(\bA)$. Hence,
$$\Lambda_\infty(\bA) =\cap_{k\ge 1}\Lambda_k(\bA)
\subseteq \cap_{k\ge 1}\cl(\Lambda_k(\bA))=W_e(\bA)
\subseteq\cap_{k\ge 1}\Omega_k(\bA)=\Omega_\infty(\bA)\,.$$

To show that $\Omega_\infty(\bA)\subseteq W_e(\bA)$.
Suppose $\ba  \in \Omega_k(\bA)$. By Remark \ref{2.3},
we may assume that $\ba=\b0$. So, $\lambda_k(\bc\cdot \bA)\ge 0$
for all $k\ge 1$ and  $\bc\in\IR^m$. For
$\bF=(F_1,\dots,F_m) \in\SH^m \cap \FH^m$, let
$K=\sum_{i=1}^m\rank\(F_i\)+1$. Then
$$ \lambda_1(\bc\cdot(\bA+\bF))\ge  \lambda_{K+1}(\bc\cdot \bA)\ge 0\
\mbox{ and }\ \lambda_1(-(\bc\cdot(\bA+\bF)))\ge
\lambda_{K+1}(-\bc\cdot \bA)\ge 0\,.$$
Therefore,   $\bc\cdot\b0=0\in  \cl(W(\bc\cdot(\bA+\bF)))=\bc \cdot
\cl(W(\bA+\bF) )$. Hence, $\bc\cdot \b0\in \bc \cdot W_e(\bA)$ for all
$\bc\in \IR^m$. By the convexity of $W_e(\bA)$, we have $\b0\in W_e(\bA)$.
\qed

\begin{example} For $n\ge 1$, let
$B_n=\[\begin{array}{rr}\frac1n&0\\ 0&-\frac1n\end{array}\]$,
$C_n=\[\begin{array}{rr}\frac1n&0\\ 0&0\end{array}\]$,
$A_1=\oplus_{n=1}^\infty B_n$,
$A_2=\oplus_{n=1}^\infty C_n$ and $\bA=(A_1,A_2)$.
Then $(0,0)\in \Omega_\infty (\bA)$ but $\Lambda_\infty(\bA)=\emptyset$.
\end{example}

\noindent{\bf Acknowledgement } Li is an honorary professor of the University of Hong Kong. His research  was partially supported by an NSF grant and
the William and Mary Pulmeri Award.

\medskip
\begin{tabular}{cc}
CHI-KWONG LI&YIU-TUNG POON\\
Department of Mathematics&  Department of Mathematics \\
The College of William and Mary&Iowa State University\\
Williamsburg, Virginia 23185&Ames, Iowa 50011\\
 USA &USA\\
 ckli@math.wm.edu&ytpoon@iastate.edu\end{tabular}


\begin{thebibliography}{WWW}

\bibitem{AS} J. Anderson and J.G. Stampfli, Compressions and commutators,
Israel J. Math. 10 (1971), 433-441.


\bibitem{C} M.D. Choi,  Completely positive linear maps on complex matrices,
Linear Algebra and Appl. 10 (1975), 285-290.

\bibitem{Cet}
M.D. Choi, M. Giesinger, J. A. Holbrook, and D.W. Kribs,
Geometry of higher-rank numerical ranges,
Linear and Multilinear Algebra 56 (2008), 53-64.


\bibitem{Cet0}
M.D. Choi, J.A. Holbrook, D. W. Kribs, and K. {\.Z}yczkowski,
Higher-rank numerical ranges of unitary and normal matrices,
Operators and Matrices 1 (2007), 409-426.

\bibitem{Cet1}
M.D. Choi, D. W. Kribs, and K. {\.Z}yczkowski,
Higher-rank numerical ranges and compression problems,
Linear Algebra Appl., 418 (2006), 828--839.

\bibitem{Cet2}
M.D. Choi, D. W. Kribs, and K. {\.Z}yczkowski,
Quantum error correcting codes from the compression formalism,
Rep. Math. Phys., 58 (2006), 77--91.


\bibitem{GLW} H.L. Gau, C.K. Li, and P.Y. Wu,
Higher-Rank Numerical Ranges and Dilations,
J. Operator Theory, to appear.

\bibitem{GR} K.E. Gustafson and D.K.M. Rao, Numerical ranges:
The field of values of linear operators and matrices, Springer,
New York, 1997.

\bibitem{KL} E. Knill and R. Laflamme,
Theory of quantum error correcting codes, Phys. Rev. A 55 (1997), 900-911.


\bibitem{KLV} E. Knill, R. Laflamme,  and L. Viola, Theory of quantum error
correction for general noise, Phys. Rev. Lett. 84 (2000), no. 11, 2525--2528.

\bibitem{LP} C.K. Li and Y.T. Poon,
Convexity of the joint numerical range, SIAM J. Matrix
Analysis Appl. 21 (1999), 668-678.

\bibitem{LP2} C.K. Li and Y.T. Poon,
The Joint Essential Numerical Range of operators:
Convexity and Related Results, submitted.
http://www.resnet.wm.edu/$\,\tilde{}\,$cklixx/jwess.pdf.



\bibitem{LPS} C.K. Li, Y.T. Poon, and N.S. Sze,
Condition for the higher rank numerical range to be non-empty,
Linear and Multilinear Algebra, to appear.
http://arxiv.org/abs/0706.1540.

\bibitem{LPS2} C.K. Li, Y.T. Poon, and N.S. Sze,
Higher rank numerical ranges and low
rank perturbations of quantum channels,
J. Math. Anal. Appl. 348 (2008), no. 2, 843--855.


\bibitem{LS} C.K. Li and N.K. Sze, Canonical forms,
higher rank numerical ranges, totally isotropic subspaces,
and matrix equations,
Proc. Amer. Math. Soc. 136 (2008), no. 9, 3013--3023.


\bibitem{R} R.A. Martinez-Avendano, Higher-rank numerical range in
infinite-dimensional Hilbert space, Operators and Matrices 2 (2008), 249-264.

\bibitem{Tv} H. Tverberg, A generalization of Radon's theorem, J. Lond. Math. Soc. 41 (1966), 123-128.

\bibitem{W1} H. Woerdeman,
The higher rank numerical range is convex,
Linear and Multilinear Algebra 56 (2008), 65-67.

\end{thebibliography}
\end{document}